\documentclass[12pt]{article}

\usepackage{amssymb}  
\usepackage{a4}  
\newtheorem{theorem}{Theorem}[section]

\newcommand{\be}{\begin{equation}}  
\newcommand{\ee}{\end{equation}}  
\newcommand{\bea}{\begin{eqnarray}}  
\newcommand{\eea}{\end{eqnarray}}  
\newcommand{\tl}{\,\triangleleft}  
\def\cocross{{>\!\!\!\triangleleft\,}}  
\def\1{{\bf 1}}  
\def\id{\mbox{id\,}}  
\newcommand{\bphi}{\mbox{\boldmath $\phi$}} 
\newcommand{\bpsi}{\mbox{\boldmath $\psi$}} 
\def\R{\mathcal{R}\,}  
\def\b#1{{\mathbb #1}}  

\begin{document}

\title{New approach to Hermitian  $q$-differential operators on
$\b{R}_q^N$}

\author{Gaetano Fiore  \footnote{Talk given at the Conference
``Noncommutative Geometry and Representation Theory in Mathematical Physics'', 
Karlstad (Sweden), July 2004.
Preprint 04-45 Dip.
Matematica e Applicazioni, Universit\`a di Napoli;  DSF/04-35 ; ESI 1529. Work
partially supported by the European Commission RTN Programme
HPRN-CT-2000-00131 and by MIUR} \\\\           \and 
        Dip. di Matematica e Applicazioni, Fac.  di Ingegneria\\  
        Universit\`a di Napoli, V. Claudio 21, 80125 Napoli \\
        and \\
        I.N.F.N., Sezione di Napoli,\\ 
        Complesso MSA, V. Cintia, 80126 Napoli 
        } 
\date{} 


\maketitle

\abstract{
We report on our recent breakthrough \cite{Fio04} in the costruction
 for $q\!>\!0$ of Hermitean and ``tractable'' differential operators out
of the $U_qso(N)$-covariant differential calculus on the 
noncommutative manifolds $\b{R}_q^N$ 
(the socalled ``quantum Euclidean spaces'').}

\section{Introduction}

Noncommutative geometry for non-compact noncommutative manifolds 
is usually much more difficult than for compact ones, especially when
trying to proceed from an algebraic to a functional-analytical treatment (see
e.g. \cite{MasNakWor03}). Major sources of complications are
$\star$-structures and $\star$-representations  of the involved algebras.
As a noncommutative space  here we adopt
the socalled $N$-dimensional quantum Euclidean space $\b{R}_q^N$ ($N\ge 3$), a
deformation of  $\b{R}^N$  characterized  by its covariance under the quantum
group $U_qso(N)$ (with $q$ real so that the latter is compact; for
$\star$-representation purposes we shall later need and choose $q>0$) and
shall denote by $F$ the $\star$-algebra ``of functions on $\b{R}_q^N$''. The
$\star$-structure \cite{Ogi92} associated to the corresponding
$U_qso(N)$-covariant differential calculus \cite{CarSchWat91}  on  $\b{R}_q^N$
\cite{FadResTak89} is however characterized   by an unpleasant  nonlinear
action  on the  differentials, the partial derivatives and the exterior
derivative.  This at the origin of a host  of formal and substantial
complications: it seems hardly possible to construct physically relevant
differential operators (e.g. momentum components, Laplacian, kinetic term in a
would-be field theory on $\b{R}_q^N$) which are at the same time Hermitean and
``tractable''.

Here we report on a recent breakthrough \cite{Fio04} to the
approach to these problems. This is based 
on our result (Thm 1 in \cite{Fio04}) that the transformation
laws \cite{Ogi92} of the $N$ partial derivatives  $\partial^{\alpha}$ and of
the exterior derivative $d$ of the differential
calculus \cite{CarSchWat91} on  $\b{R}_q^N$  under the $\star$-structure can
be expressed as the similarity transformations  
\be 
\partial^{\alpha}{}^{\star}= -\tilde \nu '{}^{-2}\tilde\partial^{\alpha} 
\tilde \nu '{}^2,  \qquad \qquad \qquad d{}^{\star}= -\tilde \nu '{}^2
d \tilde \nu'{}^{-2},                       \label{simila} 
\ee 
where $\tilde \nu '{}$ is a  $U_qso(N)$-invariant
positive-definite pseudodifferential operator, more precisely  
the realization of the fourth 
root of the ribbon element of the extension of  
$U_qso(N)$ with a central element generating  
dilatations of $\b{R}_q^N$.

The situation reminds us of what one encounters  in functional analysis 
on the real line when taking 
 the Hermitean conjugate of a differential operator like 
$$
\partial:=\tilde \nu'{}(x)\, \frac d{dx}\,\frac 1{\tilde\nu'{}(x)}, 
   \label{undef} 
$$
where $\tilde\nu'{}(x)$ is a positive  smooth function; 
as an element of the Heisenberg algebra $\partial$  
is not imaginary 
(excluding the trivial case  $\tilde\nu'{}\equiv 1$) 
 w.r.t. the $\star$-structure  
$$ 
x^{\star}=x, \qquad \left(\frac d{dx}\right)^{\star}=-\frac d{dx}, 
$$ 
but fulfills the similarity transformation (\ref{simila}).
This corresponds to the fact that it is not 
antihermitean as an operator 
on $L^2(\b{R})$. $\partial$ is however (formally) antihermitean 
on $L^2(\b{R},\tilde\nu '{}^{-2}dx)$, the Hilbert space of square-integrable
functions over  $\b{R}$  endowed with measure $\tilde\nu '{}^{-2}dx $.  
In other words, if we  introduce the scalar product
\be
(\bphi,\bpsi)=\int \bphi^{\star}(x)\,\tilde\nu '{}^{-2}(x)\,\bpsi(x)dx, 
                                                          \label{weight} 
\ee
[as one does when setting the Sturm-Liouville problem 
for  $\partial^2$], $\partial$ becomes antihermitean under the corresponding Hermitean 
conjugation $\dagger$\footnote{The Hermitean conjugation $\dagger$ is 
the representation of the following modified $\star$-structure $\star'$ of the 
Heisenberg algebra:  $a^{\star'}=(\tilde\nu '{}^{-2}\,a \tilde\nu '{}^2)^{\star}=   
\tilde\nu '{}^2\,a^* \tilde\nu '{}^{-2}$. 
}: 
$$ 
( A^{\dagger}\bphi,\bpsi):=( \bphi,A\bpsi )
\qquad\quad\Rightarrow\qquad\quad \partial^{\dagger}=-\partial. 
$$
The idea is to adopt the analog of (\ref{weight}), where now $\nu'$ is a 
{\it pseudodifferential operator} rather than a function,
as a scalar product also in our noncommutative space.

\section{Preliminaries and notation}

\subsection*{The algebras.} $F$ is essentially the unital associative 
algebra over $\b{C}[[h]]$ generated by $N$ elements $x^{\alpha}$ (the cartesian 
``coordinates'') modulo the relations  
(\ref{xxrel}) given below, and 
is extended to include formal power series in the generators; 
out of $F$ one can extract subspaces consisting of elements 
that can be considered integrable or square-integrable functions. 
The $U_qso(N)$-covariant differential calculus 
on $\b{R}_q^N$ \cite{CarSchWat91} is defined introducing  
the invariant exterior derivative $d$, satisfying nilpotency 
and the Leibniz rule $d(fg)=dfg+fdg$, and imposing the covariant 
commutation relations 
 (\ref{xxirel}) between the $x^{\alpha}$ and the differentials 
$\xi^{\alpha}:=dx^{\alpha}$.  
 Partial derivatives are introduced through the  
decomposition 
$d=:\xi^{\alpha}\partial_{\alpha}$. All the other commutation relations 
are derived by consistency. The complete list is 
given by 
\bea 
&& \mathcal{P}_a{}^{\alpha\beta}_{\gamma\delta}x^{\gamma}x^{\delta}=0, \label{xxrel}\\ 
&& x^{\gamma}\xi^{\alpha}=q\hat R^{\gamma\alpha}_{\beta\delta}\xi^{\beta}x^{\delta},\label{xxirel}\\ 
&& (\mathcal{P}_s+\mathcal{P}_t)^{\alpha\beta}_{\gamma\delta}\xi^{\gamma}\xi^{\delta}=0,\label{xixirel}\\ 
&& \mathcal{P}_a{}^{\alpha\beta}_{\gamma\delta}\partial_{\beta}\partial_{\alpha}=0, \label{ddrel}\\ 
&& \partial_{\alpha} x^{\beta} = \delta^{\beta}_{\alpha}+q\hat R^{\beta\gamma}_{\alpha\delta} 
x^{\delta}\partial_{\gamma},                    \label{dxrel}\\ 
&& \partial^{\gamma}\xi^{\alpha}=q^{-1}\hat R^{\gamma\alpha}_{\beta\delta}\xi^{\beta} 
\partial^{\delta}.\label{dxirel} 
\eea 
Here the $N^2 \times N^2$ matrix $\hat{R}$ is the braid matrix of  
$SO_q(N)$ \cite{FadResTak89}. The matrices 
$\mathcal{P}_s$, $\mathcal{P}_a$, $\mathcal{P}_t$ are $SO_q(N)$-covariant 
deformations of the symmetric trace-free, 
antisymmetric and trace projectors respectively, which appear in the 
orthogonal projector decomposition of $\hat{R}$
\be 
\hat R = q\mathcal{P}_s - q^{-1}\mathcal{P}_a + q^{1-N}\mathcal{P}_t.        
                                                \label{projectorR} 
\ee 
Thus they satisfy the equations (with $A,B = s,a,t$)
\be 
\mathcal{P}_{A}\mathcal{P}_{B} = \mathcal{P}_A
\delta_{AB}, \qquad \qquad  \mathcal{P}_s+\mathcal{P}_a+\mathcal{P}_t = \1_{N^2}.
\ee 
The $\mathcal{P}_t$ projects on a
one-dimensional sub-space and   can be written in the form  
\be 
\mathcal{P}_t{}_{\gamma\delta}^{\alpha\beta} \sim g^{\alpha\beta}
g_{\gamma\delta}  ;             \label{Pt}  
\ee  
here the $N \times N$ matrix $g_{\alpha\beta}$
is a $SO_q(N)$-isotropic tensor,  
deformation of the ordinary Euclidean metric. The elements
$$ 
r^2\equiv x \cdot x:= x^{\alpha} g_{\alpha\beta}x^{\beta}, \qquad
\qquad   \partial\cdot\partial:= 
g^{\alpha\beta}\partial_{\beta} \partial_{\alpha}
$$  
are $U_qso(N)$-invariant and respectively generate the centers 
of $F,F'$. $r^2$ is real and positive-definite w.r.t. the $\star$-structure
in $F$ for real $q$ introduced in \cite{FadResTak89}.
$\partial\cdot\partial$
is a deformation of  the Laplacian on $\b{R}^N$. 
We shall adopt here mainly a set $\{x^{\alpha}\}$ of {\it real}
generators $x^{\alpha}=V^{\alpha}_ix^i$,  
whereas by $\{x^i\}$ we mean the standard non-real basis adopted since the
works  \cite{FadResTak89,CarSchWat91}, in which $r^2$ takes the form
$r^2=\sum_i x^i{}^{\star}x^i$ and the relations (\ref{xxrel}-\ref{dxirel})
take  the simplest explicit form, and $V$ is the (complex) linear
transformation relating the two bases.\footnote{For instance, 
relations
(\ref{xxrel}), (\ref{xixirel})  for $\b{R}_q^3$ in the basis
$\{x^i\}=\{x^-,x^0,x^+\}$ \cite{FadResTak89,CarSchWat91} amount to: 
\bea
&&x^-x^0=qx^0x^-\qquad x^0x^+=qx^+x^0,\qquad 
[x^-,x^+]=(q^{\frac 12} -q^{-\frac 12} )x^0x^0 ,\nonumber\\[2pt]
&&q\xi^- \xi^0 +\xi^0 \xi^-=q\xi^0
\xi^+ +\xi^+ \xi^0=\xi^- \xi^+ + \xi^+ \xi^-=(\xi^{\pm})^2=
(\xi^0)^2 -(q^{\frac 12} -q^{-\frac 12})\xi^-\xi^+= 0.\nonumber
\eea
The $\star$-structure acts on the coordinates as follows:
$$ 
(x^-)^* =  q^{\frac 12}x^+, \qquad (x^0)^* = x^0, \qquad 
(x^+)^* = q^{-\frac 12}x^-. 
$$
}

We shall slightly extend $F$ by  introducing the square root $r$ of $r^2$  
and its inverse $r^{-1}$ as new (central) generators; 
$r$ is ``the deformed Euclidean 
distance  of the generic point of coordinates $(x^i)$ of $\b{R}_q^N$  
from the origin''. The elements 
$$ 
t^{\alpha}:=r^{-1}x^{\alpha}
$$ 
fulfill $t\cdot t=1$ and thus generate the subalgebra $F(S_q^{N\!-\!1})$ 
of ``functions on the unit quantum Euclidean sphere''.  $F(S_q^{N\!-\!1})$  can be 
completely decomposed into eigenspaces $V_l$ of the quadratic Casimir of 
$U_qso(N)$,  or equivalently of the Casimir $w$ defined in (\ref{defw}) 
with eigenvalues $w_l:=q^{-l(l+N-2)}$, implying a corresponding  
decomposition for $F$: 
\be 
F(S_q^{N\!-\!1})=\bigoplus\limits_{l=0}^{\infty}V_l, \qquad\qquad 
F=\bigoplus\limits_{l=0}^{\infty} 
\left(V_l\otimes \b{C}[[r,r^{-1}]]\right). \label{ldeco}  
\ee 
An orthonormal basis $\{S_l^I\}$ (consisting of
`spherical harmonics') of $V_l$  
can be extracted out of the set of homogeneous, 
completely symmetric and trace-free polynomials 
of degree $l$, suitably normalized. 
Therefore for the generic $f\in F$ 
\be 
f=\sum\limits_{l=0}^{\infty} f_l=\sum\limits_{l=0}^{\infty} 
\sum_IS_l^I f_{l,I}(r).                       \label{ldeco2}  
\ee 
 
We shall call $\mathcal{DC}^*$ (differential calculus algebra 
on $\b{R}_q^N$) the unital associative  
algebra  over $\b{C}[[h]]$ generated by  
$x^{\alpha},\xi^{\alpha},\partial_{\alpha}$ modulo relations
(\ref{xxrel}-\ref{dxirel}).  We shall denote by  
$\bigwedge^*$ (exterior algebra, or algebra of exterior forms) 
the graded unital subalgebra generated by the $\xi^{\alpha}$ alone, with 
grading $\natural \equiv$the degree in $\xi^{\alpha}$, and by $\bigwedge^p$  
(vector 
space of exterior $p$-forms) the component with grading $\natural =p$,  
$p=0,1,2...$. 
Each $\bigwedge^p$ carries an irreducible representation 
of $U_qso(N)$, and its dimension is the binomial 
coefficient $N\choose{p}$ \cite{Fio04jpa}, exactly as in the $q=1$ 
(i.e. undeformed) case.
We shall endow $\mathcal{DC}^*$ with the same grading $\natural $, and call 
$\mathcal{DC}^p$ its component with grading $\natural =p$. The elements of 
$\mathcal{DC}^p$ can be considered differential-operator-valued $p$-forms. 
We shall denote by $\mathcal{H}$ (Heisenberg algebra on $\b{R}_q^N$) 
the unital subalgebra generated by the $x^{\alpha}, \partial_{\alpha}$. Note that 
by definition $\mathcal{DC}^0=\mathcal{H}$, and that both 
$\mathcal{DC}^*$ and $\mathcal{DC}^p$ are $\mathcal{H}$-bimodules. 
We shall denote by $\Omega^*$ (algebra of differential forms) 
the graded unital subalgebra generated by the $\xi^{\alpha},x^{\alpha}$, with 
grading $\natural $, and by $\Omega^p$  
(space of differential $p$-forms) its component with grading $p$;  
by definition $\Omega^0=F$ itself. 
Clearly both $\Omega^*$ and $\Omega^p$ are 
$F$-bimodules. 
Finally, we recall that  by replacing $q\to q^{-1}, \hat R\to \hat R^{-1}$ in 
(\ref{xxrel}-\ref{dxirel}) one gets an alternative $U_qso(N)$-covariant
calculus, whose exterior/partial derivative and differentials  
$\hat d,i\hat\partial^{\alpha},\hat\xi^{\alpha}$ are the $\star$-conjugates
of $d,i\partial^{\alpha},\xi^{\alpha}$ and are related to the latter by a 
rather complicated rational
transformation \cite{OgiZum92}.

\subsection*{Hodge duality.} The ``Hodge map'' 
\cite{Fio04jpa} is a $U_qso(N)$-covariant,  $\mathcal{H}$-bilinear map 
\be 
*:\mathcal{DC}^p\to\mathcal{DC}^{N-p} 
\ee 
($p=0,1,...,N$) such that ${}^*\1=dV$ and on each $\mathcal{DC}^p$  (hence
on the whole $\mathcal{DC}^*$) 
\be 
*^2\equiv *\circ *=\id,       \label{involution} 
\ee 
defined by setting on the monomials in the $\xi^{\alpha}$
\be 
{}^*(\xi^{\alpha_1}...\xi^{\alpha_p})=q^{-N(p-N/2)} 
c_p\,\xi^{\alpha_{p+1}}...\xi^{\alpha_N} 
\varepsilon_{\alpha_N...\alpha_{p+1}}{}^{\alpha_1...\alpha_p}\Lambda^{2p-N}.
 \label{defHodge2} 
\ee 
Here $c_p$ is a suitable normalization factor [going to $1/(N\!-\!p)!$
in the commutative limit $q\to 1$], $\varepsilon^{\alpha_1...\alpha_N}$ is the 
$q$-deformed $\varepsilon$-tensor
and $\Lambda^{\mp 1}$ are the square
root and the inverse square root of
\be 
\Lambda^{-2}:=1 +(q^2-1)x^{\alpha}\partial_{\alpha}+ 
\frac{q^{N-2}(q^2-1)^2}{(1+q^{N-2})^2}r^2\partial\cdot\partial\equiv 1+O(h); 
\ee 
as a pseudodifferential operator $\Lambda$ acts as a dilatation:
\be 
\Lambda x^{\alpha}=q^{-1}x^{\alpha}\Lambda,\qquad\qquad 
\Lambda\partial^{\alpha}=q\partial^{\alpha}\Lambda, \qquad\qquad 
\Lambda \xi^{\alpha}=\xi^{\alpha}\Lambda,  \qquad\qquad  \Lambda 1=1.\label{Lambdaprop} 
\ee 
In $\mathcal{DC}^*$ one can introduce \cite{Fio04,CerFioMad00} as an
alternative basis of 1-forms a ``frame'' \cite{DimMad96} $\{\theta^{\alpha}\}$, in the sense 
 that $[\theta^{\alpha},\mathcal{H}]=0$
(implying in particular $[\theta^{\alpha},F]=0$). In the latter
(\ref{defHodge2} ) takes the form 
\be 
{}^*(\theta^{\alpha_1}\theta^{\alpha_2}...\theta^{\alpha_p}) 
=\,c_p\,\theta^{\alpha_{p+1}}...\theta^{\alpha_N} 
\varepsilon_{\alpha_N...\alpha_{p+1}}{}^{\alpha_1...\alpha_p},  \label{defHodge1} 
\ee 

Restricting the domain of $*$ to the unital 
subalgebra $\widetilde{\Omega}^*\subset \mathcal{DC}^*$ generated by 
$x^{\alpha},\xi^{\alpha},\Lambda^{\pm 1}$ one obtains a $U_qso(N)$-covariant,
$\widetilde{F}$-bilinear map (with $\widetilde{F}\equiv \widetilde{\Omega}^0$)
\be
*:\widetilde{\Omega}^p\to\widetilde{\Omega}^{N-p}    \label{restri}
\ee
fulfilling again ${}^*\1=dV$  and (\ref{involution}).  Finally,
introducing the exterior coderivatives
\be
\delta:= -{}^*\, d\,{}^*\qquad\qquad\hat\delta:= -{}^*\, \hat d\,{}^*
\ee 
one finds that
on all of $\mathcal{DC}^*$, and in particular on all of $\Omega^*$,
the Laplacians 
\bea
\hat\Delta &:=& d\,\delta+\delta\,d=-q^2\,\partial \cdot \partial \Lambda^2=
-q^{-N}\hat\partial \cdot \hat\partial,         \label{Laplacian}\\
\Delta &:=&\hat d\,\hat
\delta\!+\!\hat \delta\,\hat d=-q^{-2}\,\hat \partial \cdot \hat \partial \Lambda^{-2}
=-q^N\,\partial \cdot \partial
\eea
are respectively quadratic in the $\hat\partial_{\alpha},\partial_{\alpha}$ 
(note: not in the $\partial_{\alpha},\hat\partial_{\alpha}$!).

The restriction (\ref{restri}) is the notion closest to the conventional notion
of a Hodge map on $\b{R}_q^N$: as a matter of fact, there is no
$F$-bilinear restriction of $*$ to $\Omega^*$.  
However  $\widetilde{\Omega}^*$ is not closed 
under the $\star$-structure which we shall recall below.

\subsection*{$\widetilde{U_qso(N)}$ and its realization by
(pseudo)differential operators.}  
We extend \cite{Maj94} the Hopf algebra $U_qso(N)$   by adding  
a central, primitive generator $\eta$ 
$$ 
\phi(\eta)=\1\otimes\eta+\eta\otimes\1, \qquad\qquad 
\epsilon(\eta)=0,\qquad\qquad S\eta=-\eta 
$$ 
(here $\phi,\epsilon,S$ respectively denote the coproduct, counit, antipode),
and we endow the resulting Hopf algebra $\widetilde{U_qso(N)}$  
by the  quasitriangular structure
\be 
\tilde\R:=\R\, q^{\eta\otimes\eta}, 
\ee 
where $\R\equiv\R^{(1)}\otimes\R^{(2)}$ (in a Sweedler notation with upper
indices and suppressed  summation index) denotes the  quasitriangular
structure of $U_qso(N)$. The action (which here and in \cite{Fio04} we choose
to be {\it right}) on $\mathcal{DC}^*$ is completely specified by the
transformation laws of the generators
$\sigma^{\alpha}=x^{\alpha},\xi^{\alpha},\partial^{\alpha}$,   which read
\bea
&&\sigma^{\alpha}\tl \, g=\rho^{\alpha}_{\beta}(g)\sigma^{\beta}\qquad\qquad  
\sigma^{\alpha}=x^{\alpha},\xi^{\alpha},\partial^{\alpha}, 
\qquad g\in U_qso(N), 
\nonumber\\ 
&&x^{\alpha}\tl \,\eta=x^{\alpha},\qquad\qquad \xi^{\alpha}\tl \,\eta=\xi^{\alpha}, 
\qquad\qquad \partial^{\alpha}\tl \,\eta=-\partial^{\alpha}; \nonumber
\eea 
here $\rho$ denotes the $N$-dimensional representation of $U_qso(N)$. The
elements 
$$
Z^{\alpha}_{\beta}:= T^{(1)}\!\rho^{\alpha}_{\beta}(T^{(2)}), \qquad \qquad       
T:=\R_{21}\R\equiv T^{(1)}\!\otimes\! T^{(2)},  
\qquad\quad\R_{21}\equiv\R^{(2)}\!\otimes\!\R^{(1)}
$$  
are generators of $U_qso(N)$, and make up the
``$SO_q(N)$ vector field matrix''
$Z$   \cite{Zum91,SchWatZum92}.
We recall that the ribbon element $w\in U_qso(N)$ 
is the central element such that  
$$
 w^2=u_1 S(u_1), \qquad\qquad 
u_1:= (S\R^{(2)}) \R^{(1)}. \label{defw}\\ 
$$
It is well-known \cite{Dri90} that there exist isomorphisms
$U_hso(N)[[h]]\simeq Uso(N))[[h]]$ of $\star$-algebras over $\b{C}[[h]]$. This
essentially means that it is possible to express the elements of either
algebra as power series in $h=\ln q$ with coefficients in the other. In
particular $w$ has an extremely simple expression in terms of 
 the quadratic Casimir  $C$ of $so(N)$: 
\be
w=q^{-C}=e^{-hC}=\1+O(h), \qquad\qquad\qquad C:=X^aX_a=:L(L+N-2).
\ee
($\{X^a\}$ is a basis of $so(N)$).
We denote by $\nu:=w^{1/4}$ and by $\tilde w,\tilde \nu,\tilde T,\tilde Z$ 
the analogs of $w,\nu,T,Z$ obtained by replacing 
$\R$ by $\tilde\R$. As an immediate  consequence 
$$ 
\tilde w=q^{-C}q^{-\eta^2}, \qquad 
\tilde \nu= q^{-C/4}q^{-\eta^2/4},\qquad \tilde T= T q^{2\eta\otimes\eta}.  
$$ 

As shown in \cite{Fio95cmp,ChuZum95,Fio04jpa}, for real $q$ there exists a 
$\star$-algebra homomorphism\footnote{
This is the $q$-deformed analog of the realization  in the $q=1$ case of the
Cartan-Weyl $L^{\alpha\beta}$ generators of $\widetilde{Uso(N)}$ through
vector fields (i.e. $1^{st}$-order differential operators):
$$
\varphi(L^{\alpha\beta})=x^{\alpha}\partial^{\beta}-x^{\beta}
\partial^{\alpha},\qquad  \varphi(x^{\alpha})=x^{\alpha}
\qquad \varphi(\partial_{\alpha})=\partial_{\alpha}, \qquad 
\varphi(\eta)=-x^{\alpha}\partial_{\alpha}.
$$
}
\be 
\varphi: \mathcal{H}\cocross \widetilde{U_qso(N)}\to \mathcal{H},              \label{Hom1} 
\ee 
acting as the identity on $\mathcal{H}$ itself, 
\be 
\varphi(a)=a \qquad\qquad a\in \mathcal{H}.     \label{Hom1'} 
\ee 
(This requires introducing an additional generator  
$\eta'=\varphi(\eta)\in\mathcal{DC}^*$ subject to the condition 
$\varphi(q^{\eta})=q^{\eta'}=q^{-N/2}\Lambda$, so that
$[\eta',\xi^{\alpha}]=0$, $[\eta',x^{\alpha}]=-x^{\alpha}$,
 $[\eta',\partial^{\alpha}]=\partial^{\alpha}$. 
In the sequel we shall often use the the shorthand notation 
$$
\varphi(g)=:g',  \qquad \qquad \qquad g\in \widetilde{U_qso(N)}.  
$$
One finds \cite{Fio04} on the spherical harmonics of level $l$
(with $l=0,1,2,...$) 
\be
\nu'S_l^I =q^{-l(l+N-2)/4}\,S_l^I=S_l^I\tl \nu.  \label{eigenvalue}
\ee

\section*{$\star$-Structure and Hermiticity in configuration space
representation}

 \begin{theorem}\cite{Fio04} For $q>0$ 
the $\star$-structure of $\mathcal{DC}^*$
given in  \cite{FadResTak89,OgiZum92} can be 
expressed in the form 
\bea 
x^{\alpha}{}^{\star} &=& x^{\alpha},\hfill \\ 
\xi^{\alpha}{}^{\star} 
&=& q^N\,\xi^{\beta}\,Z'{}^{\alpha}_{\beta}\Lambda^{-2}, 
\label{starsimixi}\\ 
\partial^{\alpha}{}^{\star}&=&- \tilde \nu'{}^{-2}\partial^{\alpha}
\tilde \nu'{}^2=-q^{\frac {1-N}2} \tilde \nu'{}^{-2}\partial_{\alpha}
\tilde \nu'{}^2\Lambda,                        \label{starsimid} \\
d^{\star}&=& - \tilde \nu'{}^2d\tilde \nu'{}^{-2}. \label{starsimd} 
\eea 
\label{starsimi} 
\end{theorem} 
 
\subsection*{Integration over $\b{R}_q^N$.} 
Up to normalization this is defined by \cite{Ste96}
\be
\int_qf(x)\,d^N\!x =\int\limits^{\infty}_0dr\, m(r)\, 
r^{N\!-\!1} \int_{S_q^{N\!-\!1}}d^{N\!-\!1}\!t\,f(t,r)\nonumber\\ 
=\int\limits^{\infty}_0dr\,   m(r)\,r^{N\!-\!1} f_0(r). \nonumber 
\ee 
where $f_0$ denotes the $l=0$ component in (\ref{ldeco2}).
The`weight' $m$ along the radial direction $r$ has to
fulfill, beside $m(r)>0$, the $q$-scaling condition
\be 
m(r)=m(qr). 
\ee
Among the weights fulfilling the latter we can single out two ``extreme''
cases: 
\begin{enumerate}
\item $m(r)\equiv 1$; 
 
\item  $m(r)\equiv m_{J,r_0}(r):=|q-1|\sum\limits_{n=-\infty}^{\infty} 
r\delta(r-r_0q^n)$ (`Jackson' integral). 
\end{enumerate}
Integration over  $\b{R}_q^N$ fulfills the standard properties
of reality, positivity, $\widetilde{U_qso(N)}$-invariance, 
and a slightly modified cyclicity property \cite{Ste96}.
Moreover,  if $f$ is a regular function decreasing
faster than  $1/r^{N-1}$ as  $r\to\infty$ the Stokes  theorem holds 
\be 
\int_q \partial_{\alpha} f(x)\,d^N\!x=0, 
\qquad\qquad \int_q \hat\partial_{\alpha} f(x)\,d^N\!x=0. \label{stokes} 
\ee 
Integration of functions immediately leads to integration 
of $N$-forms $\omega_N$: upon moving all the $\xi$'s to the right 
of the $x$'s  we can express  
$\omega_N$ in the form $\omega_N=fd^N\!x$ ($d^N\!x$ denotes the unique
independent exterior $N$-form), and just
have to set   
\be 
\int_q \omega_N=\int_q f\,d^N\!x. 
\ee 

We introduce the scalar product of  two ``wave-functions'' 
$\bphi,\bpsi\in F$ and more generally  of  two
``wave-forms'' $\alpha_p, \beta_p\in \Omega^p$ by  
\be
(\bphi ,\bpsi) := \int_q \bphi^{\star} \tilde \nu '{}^{-2}
\bpsi\,d^N\!x,  \qquad\qquad                            
(\alpha_p, \beta_p):=\int_q \alpha_p^{\star}\:\:{}^*\:\tilde \nu
'{}^{-2}\beta_p.                \label{scalprod}
\ee
Setting $U^{-1}{}^{\delta}_{\gamma}:=g^{\lambda\delta}g_{\lambda\gamma}$,
the second can be easily expressed in terms of integrals 
\be
(\alpha_p, \beta_p)  
=\frac{1}{c_{N\!-\!p}}\int_q \alpha^{\star}_{\gamma_p...\gamma_1}
U^{-1}{}^{\delta_p}_{\gamma_p}...U^{-1}{}^{\delta_1}_{\gamma_1} 
\tilde \nu'{}^{-2}\beta_{\delta_p...\delta_1}{}\,d^N\!x  \label{scalprodp}
\ee
involving the form components of the generic $\omega_p\in\Omega^p$ in the
frame basis:  
$$
\omega_p =\theta^{\gamma_1}...\theta^{\gamma_p}\omega_{\gamma_p...\gamma_1}.
$$
[Note that for $q\neq 1$ $\omega_{\gamma_p...\gamma_1}\in\mathcal{H}\setminus
F$, because  $\theta^a\in\mathcal{DC}^*\setminus\Omega^*$.
As an integrand funtion at the rhs(\ref{scalprodp}) we mean what one
obtains after letting all $\partial_{\alpha}$ present under the integral
sign act  on functions on their right, using the derivation rule
(\ref{dxrel})].

As a consequence of Theorem \ref{starsimi} we obtain 
\be
(\alpha_p, \hat d\beta_{p\!-\!1}) = (\hat \delta\alpha_p, \beta_{p\!-\!1}),  
\qquad\qquad
(\hat d\beta_{p\!-\!1},\alpha_p) = (\beta_{p\!-\!1},\hat \delta\alpha_p).
\ee
and the (formal) hermiticity of both the
momenta $p^{\alpha}=i\partial^{\alpha}$ and the Laplacian
$\Delta$:
\be
(\bphi ,p^{\alpha}\bpsi) =   (p^{\alpha}\bphi ,\bpsi),   \qquad \qquad    
(\alpha_p,\Delta \beta_p) =(\Delta \alpha_p,\beta_p). \label{hermi} 
\ee
The kinetic term in the action for a (real) $p$-form (i.e. an antisymmetric 
tensor with $p$-indices) field theory of mass $M$ can be now introduced as
\be
\mathcal{S}_k=(\Delta \alpha_p,\alpha_p)
+ M^2 (\alpha_p,\alpha_p).
\ee
Spectral analysis for $p^{\alpha}, \Delta$ is affordable (see e.g.
\cite{Fio95JMP}), because of the relatively simple derivation
rules (\ref{dxrel}), compared e.g. to those of 
$p^{\alpha}_R:=p^{\alpha}\!+\!p^{\alpha}{}^{\star}$.

\medskip
The above  considerations are still {\it formal} unless we make sense out of 
$$
\tilde \nu'{}=q^{-C'/4}q^{-\eta'{}^2/4}
$$  
and its inverse as {\it positive-definite}
(pseudodifferential) operators on suitable
Hilbert subspaces  of $F,\Omega^*$. Now
$q^{-C'/4}$ is well-defined and positive-definite
by (\ref{ldeco}), (\ref{eigenvalue}).
To  define the action of  
$q^{-\eta'{}^2/4}$ on functions $\phi(r)$
 we  change  variables  
$r\to y=\ln r$, whereby $\eta'=\! -\!\partial_y\!- \!N/2$ 
and $r^{N\!-\!1}dr=e^{Ny}dy$, for any function
$\phi(r)$ denote $\tilde \phi(y):=\phi(r)$,  
and express   $e^{yN/2}\tilde 
\phi(y)=r^{N/2}\phi(r)$  in terms of its Fourier transform,
\be
e^{\frac N2 y}\tilde \phi(y)=\frac 1  
{\sqrt{2 \pi}}\int^{\infty}_{-\infty} \hat \phi(\omega) e^{i\omega  
y}d\omega.                         \label{inversefourier} 
\ee  
One finds that $q^{-\eta'{}^2/4}$  acts as multiplication by $q^{\omega^2/4}$  
on the Fourier transform:
$$
e^{\frac N2 y}q^{-\eta'{}^2/4}\tilde \phi(y)=q^{-\partial_y^2/4}e^{\frac N2 
y}\tilde \phi(y)|=\frac 1  {\sqrt{2\pi}}\int^{\infty}_{-\infty}d\omega 
\hat \phi(\omega)   e^{i\omega y} q^{\omega^2/4}.
$$
For $q\!<\!1$ the UV convergence of scalar products  is drastically
improved by $q^{\omega^2/4}$.
Introducing tentatively the `Hilbert space of square integrable
functions  on $\b{R}_q^N$'
\be
L_2^{m}:=\Big\{{\bf f}(x)\equiv
\sum_{l=0}^{\infty}\sum_I S_l^I\,f_{l,I}(r)
\,\in F \:\: |\:\: ({\bf f},{\bf f})<\infty\Big\}
\label{defL_2m}
\ee 
one finds \cite{Fio04} that all works automatically with $L_2^{m\equiv 1}$
(i.e. if we choose $m(r)\!\equiv \!1$). 
Otherwise, we have to 

\noindent
1. further require that $m$ is invariant under inversion,
$m(r^{-1})=m(r)$, and fulfills  condition (A.26) in \cite{Fio04}
in order that the scalar product $(\:,\:)$ be positive-definite;

\noindent
2. at rhs(\ref{defL_2m}) replace $F$ with the subspace $E^{m,\beta,n}$
whose elements $f_{l,I}(r)$'s
admit each an analytic continuation in the complex $r$-plane 
with poles possibly only in
\be
r_{j,k}:=q^{j\!+\!\beta}e^{i\frac{\pi(2k\!+\!1)}n}, \label{poles4}
\ee
where $\beta\in\{0,1/2\}$,  $n$ is an integer
submultiple of $N$, $k=0,1,...,n\!-\!1$,  and $j\in\b{Z}$, in order that
the hermiticity (\ref{hermi}) is actually implemented.

Condition 1. is fulfilled by a large class of weights, including
the `Jackson' ones $m_{J,q^{\beta}}$.
Condition 2. selects interesting subclasses of functions,
including $q$-special functions with ``quantized parameters'', among
the simplest one e.g.
$1/[1\!+\!(q^jr)^n]$.

\bibliographystyle{amsalpha}

\end{document}